\PassOptionsToPackage{unicode}{hyperref}
\PassOptionsToPackage{hyphens}{url}
\documentclass[]{article}
\usepackage{amsmath,amssymb,amsthm,graphicx,enumitem}
\usepackage{lmodern}
\usepackage{stmaryrd}
\usepackage{iftex}
\ifPDFTeX
  \usepackage[T1]{fontenc}
  \usepackage[utf8]{inputenc}
  \usepackage{textcomp} 
\else 
  \usepackage{unicode-math}
  \defaultfontfeatures{Scale=MatchLowercase}
  \defaultfontfeatures[\rmfamily]{Ligatures=TeX,Scale=1}
\fi
\IfFileExists{upquote.sty}{\usepackage{upquote}}{}
\IfFileExists{microtype.sty}{
  \usepackage[]{microtype}
  \UseMicrotypeSet[protrusion]{basicmath} 
}{}
\makeatletter
\@ifundefined{KOMAClassName}{
  \IfFileExists{parskip.sty}{%
    \usepackage{parskip}
  }{
    \setlength{\parindent}{0pt}
    \setlength{\parskip}{6pt plus 2pt minus 1pt}}
}{
  \KOMAoptions{parskip=half}}
\makeatother
\usepackage{setspace}
\usepackage{xcolor}
\usepackage{longtable,booktabs,array}
\usepackage{calc} 
\usepackage{etoolbox}
\makeatletter
\patchcmd\longtable{\par}{\if@noskipsec\mbox{}\fi\par}{}{}
\makeatother
\IfFileExists{footnotehyper.sty}{\usepackage{footnotehyper}}{\usepackage{footnote}}
\makesavenoteenv{longtable}
\usepackage[normalem]{ulem}
\ifLuaTeX
  \usepackage{selnolig}  
\fi
\IfFileExists{bookmark.sty}{\usepackage{bookmark}}{\usepackage{hyperref}}
\IfFileExists{xurl.sty}{\usepackage{xurl}}{} 
\hypersetup{
  hidelinks,
  pdfcreator={LaTeX via pandoc}}
\author{}
\date{}

\begin{document}
\singlespacing  
\begin{center}
\section{\textbf{An Anomaly in Diagonalization}}
T. Parent (Nazarbayev University)\\
nontology@gmail.com

\end{center}

\textbf{0. Introduction}

Primitive recursive (p.r.) functions are among the best-understood classes of computable functions, and it has long been established that no universal p.r. function exists. The standard impossibility proof relies on diagonalization to show that such a function would generate contradiction. However, it is known that some diagonal arguments, when formalized, do not demonstrate the impossibility of the diagonal object, but instead reveal a breakdown in definability or encoding. For example, in a formal setting, Richard’s paradox does not yield a contradiction; it instead reflects that one of the relevant sets is ill-defined. (For elaboration and further examples, see Simmons 1993, Chapter 2.)

This invites the possibility that other diagonal arguments may reflect similar anomalies. The argument against a universal p.r. function is considered in this light. The impetus is an algorithm which appears to satisfy all standard criteria for being p.r. while simulating the computation of $f_i(i,n)$ for any index $i$ of a binary p.r. function. The paper does not attempt to explain why this construction apparently survives the usual diagonal objection, but presents it in a form precise enough to support that analysis.

One preliminary: The following normally does not discuss p.r.
functions directly but rather function symbols
\(f_{\underline{m}}^{\underline{n}}\) within a version of Primitive
Recursive Arithmetic (PRA).\footnote{This makes the paper reminiscent of work by Princeton mathematician Edward Nelson (2011), in which PRA was purportedly shown inconsistent. (Tao's 2011 response caused Nelson 2011a to retract his claim.) This paper is emphatically not alleging that PRA is inconsistent. Rather, the claim is that the significance of one diagonal argument has been misidentified. (A more minor difference is that Nelson was focused on the induction axiom, whereas this paper concerns more minimal versions of PRA which do not include the induction axiom.)}
This is because the algorithm relies crucially on
numeric subscripts, and it is unproblematic which numeric subscript a
function symbol has. In contrast, it is dubious to assume that,
e.g., the factorial function itself has some specific index
attached to it (in Plato's heaven, as it were). The discussion thus
proceeds primarily by considering function symbols rather than functions per se.

\bigskip
\textbf{1. Exposition of PRA--}

We shall focus on a weakened version of PRA; cf. Skolem (1923); Hilbert
\& Bernays (1934, ch. 7); Curry (1941). Call it PRA--. The weakening is
that the induction axiom is replaced with the axiom that any positive
integer has a predecessor (cf. Q, Robinson's 1950
arithmetic).\footnote{Also, PRA sometimes includes principles for
  course-of-values recursion. Such things are left aside also.} But for
our purposes, the most important feature of PRA-- is the way a p.r.
function is assigned a subscripted function symbol; this is detailed
later in this section.

First, terms of PRA-- are defined inductively as follows:

\begin{enumerate}[leftmargin=1.5cm]
\def\labelenumi{(\alph{enumi})}

\item
  `0' is a term.

\item
If $\tau$ is a term, then so is $\tau^\prime$. (`0' followed by zero or more
occurrences of `$^\prime$' are the numerals. Let \uline{$\tau$} be the numeral
  co-referring with $\tau$.)\footnote{In what follows, I often elide the
    distinction between use and mention. Quine's (1951) corner quotes could be used to distinguish an expression
    $\ulcorner\tau^{\prime}\urcorner$ from what it represents.
    (Some other notation for Gödel numbers would then be required.) But
    to avoid clutter, I instead rely on context to disambiguate.}

\item
  
  \emph{v\uline{n}} is a term (a variable).
  
\item
  
If $\tau$\textsubscript{1}, \ldots, $\tau$\emph{\textsubscript{m}} are terms,
  then so is
  \emph{f\uline{nm}}($\tau$\textsubscript{1}\ldots $\tau$\emph{\textsubscript{m}})
  (a non-numeric or nn-term).
  
\end{enumerate}

But for convenience, Arabic numerals will be used and I revert to
using `\emph{x},', `\emph{y}', etc., as variables. Also, a nn-term will
normally be written as
\emph{f\textsubscript{\uline{n}}}($\tau$\textsubscript{1}, \ldots,
$\tau$\emph{\textsubscript{m}}). 

The well-formed formulae (wff) of PRA-- are defined thus:

\begin{enumerate}[leftmargin=1.5cm]
\def\labelenumi{(\roman{enumi})}

\item
  If $\tau$\textsubscript{1} and $\tau$\textsubscript{2} are terms, then
  $\tau$\textsubscript{1}=$\tau$\textsubscript{2}~is a wff.

\item

  If $\phi$ is a wff, then so is~ $\sim\!\phi$.

\item

  If $\phi$ and $\psi$ are wff, then so is ($\phi\supset\psi$).

\end{enumerate}

Assume that `=', `$\sim$', and `$\supset$', have their standard
interpretations. Parentheses will be omitted when there is no danger of
confusion.

The system has all axioms from the Hilbert schemes:

\begin{quote}
(H1) $\phi\supset(\phi\supset\phi)$

(H2) $(\phi\supset(\psi\supset\chi))\supset((\phi\supset\psi)\supset(\phi\supset\chi))$

(H3) $(\sim\psi\supset~\sim\phi)\supset(\phi\supset\psi)$

\end{quote}
PRA-- also has the standard axiomatic analysis of `=', as per the Law of
Universal Identity and the Indiscernability of Identicals:

\begin{quote}
(U\textsubscript{=}) \emph{x}=\emph{x}

(I\textsubscript{=}) $\tau\textsubscript{1}\!\!=\!\!\tau\textsubscript{2}\supset
(\phi{[}\emph{x}/\tau\textsubscript{1}{]}\supset
\phi{[}\emph{x}/\tau\textsubscript{2}{]})$
\end{quote}

Here and elsewhere, $\phi${[}\emph{x}/$\tau${]} is the result of uniformly
replacing `\emph{x}' in $\phi$ with $\tau$.

Further, in PRA-- each p.r. function is assigned to at least one basic
function symbol \(f_{\underline{n}}\); axioms defining those symbols
appropriately are included. Per Gödel (1931, pp. 179-180), the axioms
shall be specified inductively.\footnote{All references to Gödel will be
  in relation to his (1931).}

The inductive definition provided is more involved than is usual, but
the index assigned to a p.r. function will effectively code the composition of the function, and this will be crucial later. (I assume some standard way of Gödel numbering symbols, string of symbols, and sequences of
strings; cf. Gödel, pp. 179.) Moreover, it is important that the coding,
and thus, the enumeration of the axioms for the p.r. functions is
primitive recursive.\footnote{Since a p.r. enumeration must be total, let \emph{n} map to 0 when \emph{n} enumerates no axiom.} Where `*' indicates the concatenation of symbols:

\smallskip
\begin{quote}
\uline{Axioms for base p.r. functions}:

The 0\textsuperscript{th} axiom is: \emph{f}\textsubscript{0}(\emph{x})=0.
\\
The 1\textsuperscript{st} axiom is: \emph{f}\textsubscript{1}(\emph{x})=\emph{x}$^\prime$
\\
If \emph{c} codes \emph{\uline{j}}*\emph{\uline{k}} where
$1\!\leq\!\emph{j}\!\leq\!\emph{k}\!\leq\!\emph{c}$, then the
\emph{c}\textsuperscript{th} axiom is:\\
\phantom{XX}\emph{f\textsubscript{c}}(\emph{x}\textsubscript{1}, \ldots,
\emph{x\textsubscript{k}})=\emph{x\textsubscript{j}} 

\uline{Axioms for composed p.r. functions}:

If \emph{c} codes the string
\uline{0}*\emph{\uline{b}}*\uline{\emph{c}\textsubscript{1}}*
\ldots*\emph{\uline{c\textsubscript{n}}}*\emph{\uline{k}} (and each
of \emph{b}, \emph{c}\textsubscript{1}, \ldots,
\emph{c\textsubscript{n}}, and \emph{k} is less than \emph{c}), and the
\emph{b}\textsuperscript{th},
\emph{c}\textsubscript{1}\textsuperscript{th}, \ldots, and
\emph{c\textsubscript{n}}\textsuperscript{th} axioms define,
respectively,
\emph{f\textsubscript{\uline{b}}}(\emph{x}\textsubscript{1}, \ldots,
\emph{x\textsubscript{n}}),
\emph{f\textsubscript{\uline{c}}}\textsubscript{\uline{1}}(\emph{x}\textsubscript{1},
\ldots, \emph{x\textsubscript{k}}), \ldots, and
\emph{f\textsubscript{\uline{cn}}}(\emph{x}\textsubscript{1}, \ldots,
\emph{x\textsubscript{k}}), then the \emph{c}\textsuperscript{th} axiom is:\\
\phantom{XX}\emph{f\textsubscript{c}}(\emph{x}\textsubscript{1}, \ldots,
\emph{x\textsubscript{k}})=\emph{f\textsubscript{b}}(\emph{f\textsubscript{c}}\textsubscript{1}(\emph{x}\textsubscript{1},
\ldots, \emph{x\textsubscript{k}}), \ldots,
\emph{f\textsubscript{cn}}(\emph{x}\textsubscript{1}, \ldots,
\emph{x\textsubscript{k}}))

\uline{Axioms for p.r. recursions}:

If \emph{c} codes the string
`$\prime$'*\emph{\uline{a}}*\emph{\uline{d}}*\emph{\uline{k}} (and each of
\emph{a}, \emph{d}, and \emph{k} is less than \emph{c}), and the
\emph{a}\textsuperscript{th} and \emph{d}\textsuperscript{th} axioms
define, respectively,
\emph{f\textsubscript{\uline{a}}}(\emph{x}\textsubscript{1}, \ldots,
\emph{x\textsubscript{k}}) and
\emph{f\textsubscript{\uline{d}}}(\emph{x}\textsubscript{1}, \ldots,
\emph{x\textsubscript{k}}\textsubscript{+2}), then where ($\phi~\&~\psi$) is shorthand for $\sim\!(\phi~\supset~\sim\psi$), the \emph{c}\textsuperscript{th} axiom is:\\
\phantom{XX}\emph{f\textsubscript{c}}(\emph{x}\textsubscript{1}, \ldots,
\emph{x\textsubscript{k}},
0)=\emph{f\textsubscript{a}}(\emph{x}\textsubscript{1}, \ldots,
\emph{x\textsubscript{k}}) \&\\
\phantom{XX}\emph{f\textsubscript{c}}(\emph{x}\textsubscript{1}, \ldots,
\emph{x\textsubscript{k}},
\emph{y}$^\prime$)=\emph{f\textsubscript{d}}(\emph{f\textsubscript{c}}(\emph{x}\textsubscript{1},
\ldots, \emph{x\textsubscript{k}}, \emph{y}), \emph{x}\textsubscript{1},
\ldots, \emph{x\textsubscript{k}}, \emph{y})
\end{quote}

\smallskip
In both occurrences, the parenthetical ``and each\ldots is less
than \emph{c}'' is eliminable. The index \emph{c} is independently
guaranteed to exceed the arity \emph{k}, as well as the index for any
function used in defining \emph{f\textsubscript{\uline{c}}}. At least, using a standard Gödel numbering, the code for \emph{\uline{n}} is greater than \emph{n}; hence, since \emph{c} codes a string consisting of \emph{\uline{k}} and the subscripts of the component function symbols, \emph{inter alia}, \emph{c} will be greater than any of these.\footnote{Relatedly, we know that with the composed p.r. functions, if \emph{f\textsubscript{\uline{b}}} has arity \emph{n}, it must be that \emph{n}\textless{}\emph{c}. This is because \emph{c} codes \emph{n}+2 symbols, and the code for \emph{n}+2 symbols is always greater than \emph{n} itself.} This fact about \emph{c} will be important later.

Moving on, PRA-- also has the following arithmetical axioms:

\begin{quote}
(A1) $\sim\!0\!=\!\emph{x}^\prime$

(A2) $\emph{x}^\prime\!=\!\emph{y}^\prime\supset\emph{x}\!=\!\emph{y}$

(A3) $\sim\!0\!=\!\emph{x}\supset~
(\exists\emph{y}\leq\emph{x})\emph{x}\!=\!\emph{y}^\prime$

\end{quote}

As usual, `($\exists\emph{y}\leq\emph{x}$)' expresses bounded existential
quantification (a p.r. function).

\vspace{1cm}
The rules of inference in PRA-- are \emph{modus ponens} and variable
substitution:

\begin{quote}
(MP) From $\phi$ and $(\phi\supset\psi)$, $\psi$ is derivable.

(VS) From $\phi$(\emph{x}\textsubscript{1}, \ldots,
\emph{x\textsubscript{n}}),  
$\phi${[}\emph{x}/$\tau$\textsubscript{1}{]}, \ldots,
{[}\emph{x\textsubscript{n}}/$\tau$\emph{\textsubscript{n}}{]} is derivable.

\end{quote}
A finite sequence of wff counts as a \emph{derivation} of $\phi$ in PRA--
from a (possibly empty) set of wff $\Gamma$ iff: The first members are the
members of $\Gamma$ (if any), the last member is $\phi$, and any member of the
sequence is either a member of $\Gamma$, an axiom, or is derivable from
previous members via some inference rule in the system. When $\Gamma$ is empty,
we say that the sequence is a \emph{proof} of $\phi$ in PRA--.
\bigskip

\textbf{2. Schematics for a universal function}

This section begins the construction of the following function \emph{u}:

\begin{quote}
\hspace{1cm}\emph{u}(\emph{i}, \emph{n})=\emph{m} iff
\emph{f\textsubscript{\uline{i}}}(\emph{\uline{i}},
\emph{\uline{n}})=\emph{\uline{m}} is a theorem of PRA--.
\end{quote}

The existence of a universal function should be uncontroversial. But as will be argued in section 6, our construction apparently reveals that the function is p.r.

The construction starts with a p.r. proof predicate---but like
PRA, PRA-- has no predicate that represents \emph{in toto} the relations
between theorems and their proofs, as it lacks unbounded
quantifiers. If $\phi$ has Gödel number
$\ulcorner\phi\urcorner$, we can yet define an
``\emph{i}-bounded'' p.r. proof predicate \emph{B}(\emph{i}, \emph{n},
$\ulcorner\phi\urcorner$), which is satisfied
iff $\phi$ has a proof coded by \emph{n} where any axiom may occur except
those beyond the \emph{i}\textsuperscript{th} axiom defining a p.r.
function symbol.

In defining the \emph{i}-bounded proof predicate, we skip several details since they are somewhat tedious and are
identical to those for Gödel's predicate \emph{B}(\emph{x},
\emph{y}) (p. 186, \#45 in Gödel's list of p.r. functions and
relations). Or at least, only minor revisions are required. (For instance, unlike Gödel's system, PRA-- is first-order
only and has no quantifiers; many of Gödel's definitions are thus more
complex than we need.) The only revisions which
may not be obvious concern the definition of an \emph{i}-bounded
``axiom predicate'' \emph{Ax}(\emph{i}, \emph{n}) for PRA-- (cf. \#42 in
Gödel's list). This is a predicate which represents the axioms for PRA--
except those beyond the \emph{i}\textsuperscript{th} axiomatic
definition of a p.r. function symbol.

In defining \emph{Ax}(\emph{i}, \emph{n}), we first construct a p.r.
predicate \emph{AxPR}(\emph{i}, \emph{j}) which is satisfied iff
\emph{j} codes the \emph{i}\textsuperscript{th} or earlier axiom for a
p.r. function. Recall that we gave a p.r. enumeration of (the codes of) these
axioms. Hence, we know that
there is a p.r. enumeration \emph{e} such that if $\emph{e}(\emph{i})\neq0$,
\emph{e}(\emph{i}) codes the axiom defining the
\emph{i}\textsuperscript{th} p.r. function symbol.\footnote{This does
  not imply the existence of a p.r. function \emph{g}(\emph{i}) that
  outputs the \emph{i}\textsuperscript{th} p.r. function. (That alone
  would allow a contradictory diagonal function.) Our function
  does not have functions in its range, but rather codes for various
  syntactic strings. These happen to code definitions of p.r. functions.
  But a malignant diagonal function would need to \emph{compute} the
  \emph{i}\textsuperscript{th} p.r. function on input \emph{i}. And that
  would require an unbounded search of, e.g., the proofs of PRA--, to
  find one ending with an equality that identifies what
  \emph{f\textsubscript{i}}(\emph{i}) is.} This allows us to interpret the predicate as: \emph{AxPR}(\emph{i}, \emph{j}) iff
$(\exists\emph{x}\leq\emph{i})(\emph{e}(\emph{x})\!=\!\emph{j}~ \&~ \emph{j}\!\neq\!0)$.

The predicates representing the logical, arithmetical, and
`=' axioms are presumed given. So we have what we need to define
\emph{Ax}(\emph{i}, \emph{n}), and thus, to define an \emph{i}-bounded
p.r. proof predicate \emph{B}(\emph{i}, \emph{n},
$\ulcorner\phi\urcorner$) in the manner of Gödel.

Given \emph{B}(\emph{i},
\emph{n}, $\ulcorner\phi\urcorner$), we now define \emph{U}(\emph{i, n}, \emph{m}) as a predicate
that strongly represents the function \emph{u}. Briefly:

\begin{quote}
\hspace{18mm}$\emph{U}(\emph{i}, \emph{n}, \emph{m})$ iff $\emph{B}(\emph{i},
\emph{s}(\emph{i}, \emph{n}),
\ulcorner\emph{f\textsubscript{\uline{i}}}(\emph{\uline{i}},
\emph{\uline{n}})\!=\!\emph{\uline{m}}\urcorner)$
\end{quote}

This indicates  that \emph{U}(\emph{i}, \emph{n}, \emph{m}) holds iff \emph{s}(\emph{i}, \emph{n}) codes an \emph{i}-bounded proof in PRA-- of
\emph{f\textsubscript{\uline{i}}}(\emph{\uline{i}},
\emph{\uline{n}})=\emph{\uline{m}}. The suggestion will be that, if this formula has a
proof at all, the object coded by \emph{s}(\emph{i}, \emph{n})  will count as one (and otherwise it will
not).

Roughly, the idea is to define \emph{s}(\emph{i}, \emph{n}) as
follows, where \emph{e}(\emph{i}) again enumerates
(the codes of) the axioms for the function symbols:
\vspace{-12pt}
\begin{quote}
\[
s(i, n) =
\begin{cases}
    cp(\ulcorner f_{\underline{i}}(\underline{i}, \underline{n}) = \underline{m} \urcorner) & \text{if } e(i) \text{ codes an axiom defining } f_{\underline{i}}^2, \\
    0 & \text{otherwise}.
\end{cases}
\]
\end{quote}

This is meant to express that, when the defining condition is satisfied, \emph{s} outputs the code for a ``canonical
proof'' of \emph{f\textsubscript{\uline{i}}}(\emph{\uline{i}},
\emph{\uline{n}})=\emph{\uline{m}}, for some \emph{\uline{m}}.

Since \emph{f\textsubscript{\uline{i}}}(\emph{\uline{i}},
\emph{\uline{n}}) is a p.r. term, there is a unique \emph{\uline{m}} of this sort (assuming PRA--
is sound). Thus, \emph{s} effectively computes what \emph{m} is, given only
\emph{i} and \emph{n}, by generating a 
proof wherein the term
\emph{f\textsubscript{\uline{i}}}(\emph{\uline{i}}, \emph{\uline{n}}) is
reduced to a numeral. 

Our concern, however, is whether such a function can be primitive recursive. Again, diagonal arguments suggest a negative answer. Yet our algorithm for constructing the canonical proofs will appear to avoid unbounded minimization. Here is a  brief, intuitive description. (The detailed description starts in section 4.) First, the algorithm checks whether \emph{e}(\emph{i}) codes an axiom defining \(f_{\underline{i}}^{\underline{2}}\). If not, the algorithm outputs 0. Otherwise, it starts a proof with this axiom, and then instantiates it on \emph{i} and \emph{n} to produce a line that is roughly of the form
\emph{f\textsubscript{\uline{i}}}(\emph{\uline{i}},
\emph{\uline{n}})=$\phi$\emph{\textsubscript{i}}(\emph{\uline{i}},
\emph{\uline{n}}) (unless \emph{f\textsubscript{\uline{i}}} expresses a
basic binary projection, which the algorithm can handle separately).
Since the right-hand term is a p.r. term, it may then generate a proof computing the term by cycling through the primitive recursive procedures defining it in well-defined order.

But at a glance, the algorithm seems to require an unbounded search, for it is unknown how many computations are needed to reduce a given term. However, our algorithm can avoid this by utilizing a dynamically updated counter that estimates the number of computations remaining at a given stage of the proof. The counter does so by tracking the number of `\emph{f}'s that occur in each line added. The counter will often be incorrect in its estimation, but it has one redeeming feature:
It estimates the number of remaining computations to be 0 if, and only
if, the number of remaining computations \emph{is} 0. This allows the process to halt exactly when it should. Plus, the subscripts on function symbols will guide the algorithm on which p.r. operations should be applied at which stages. The result is an algorithm that seems able to compute any binary p.r. term at the requisite (primitive) recursive depth.

If so, then since a p.r. function symbol is here defined only
by function symbols with lower subscripts, a proof of
\emph{f\textsubscript{\uline{i}}}(\emph{\uline{i}},
\emph{\uline{n}})=\emph{\uline{m}} yielded by the algorithm would count
as an ``\emph{i}-bounded'' proof. In which case, \emph{s}(\emph{i},
\emph{n}) would render true the sentence $\emph{B}(\emph{i},
\emph{s}(\emph{i}, \emph{n}),
\ulcorner\emph{f\textsubscript{\uline{i}}}(\emph{\uline{i}},
\emph{\uline{n}})\!=\!\emph{\uline{m}}\urcorner)$, as intended.

\bigskip
\textbf{3. Preparatory remarks, re: function symbols within proofs}

Before detailing this, we must note well the distinctive
feature of the function symbols, namely, that the
subscripts indicate the composition of the definition for the symbol. Again, this will be crucial to the operation of the algorithm.

It may be objected that a well-defined algorithm should not be able to
``access'' the subscript on a function symbol. Officially, however, a
nn-term has no subscript; it is of the form
\emph{f\uline{nm}}($\tau$\textsubscript{1}\ldots $\tau$\emph{\textsubscript{m}}).
More importantly, the objection is at odds with the use of
\emph{AxPR}(\emph{i, n}) or Kleene's (1952) T-predicate, which are
defined with reference to indices on functions.

As an illustration of the information coded by the indices, consider the
term \emph{f\textsubscript{\uline{c}}}(\uline{1}, \uline{2}), where \emph{c} codes the
string `$\prime$'*\emph{\uline{q}}*\emph{\uline{d}}*\uline{1}. The
1\textsuperscript{st} member of the string tells us that
\emph{f\textsubscript{\uline{c}}} expresses a p.r. recursion, whereas
the 2\textsuperscript{nd} and 3\textsuperscript{rd} members indicate
that \emph{f\textsubscript{\uline{c}}} is defined by
\emph{f\textsubscript{\uline{q}}} and \emph{f\textsubscript{\uline{d}}}.
So we can recover that the axiom for \emph{f\textsubscript{\uline{c}}}
is\\ \emph{f\textsubscript{\uline{c}}}(\emph{x},
\uline{0})=\emph{f\textsubscript{\uline{q}}}(\emph{x}) \&
\emph{f\textsubscript{\uline{c}}}(\emph{x},
\emph{y}$^\prime$)=\emph{f\textsubscript{\uline{d}}}(\emph{f\textsubscript{\uline{c}}}(\emph{x},
\emph{y}), \emph{x}, \emph{y}).

Suppose also that \emph{q} codes \uline{1}*\uline{1}, whence it expresses simple identity. And suppose
\emph{d} codes
\uline{0}*\emph{\uline{p}}*\uline{1}*\emph{\uline{q}}*\emph{\uline{q}}*\uline{3}.
Since the string begins with \uline{0}, it indicates that \emph{d}
indexes a composed p.r. function---further, that
\emph{f\textsubscript{\uline{d}}} is defined by
\emph{f\textsubscript{\uline{p}}}, \emph{f}\textsubscript{\uline{1}},
and two occurrences of \emph{f\textsubscript{\uline{q}}}. So we can
recover that its axiom is \emph{f\textsubscript{\uline{d}}}(\emph{z}, \emph{x}, \emph{y})=\emph{f\textsubscript{\uline{p}}}(\emph{f}\textsubscript{1}(\emph{z}),
\emph{f\textsubscript{\uline{q}}}(\emph{x}),
\emph{f\textsubscript{\uline{q}}}(\emph{y})). (It may seem unnecessary for the axiom to include
\emph{f\textsubscript{\uline{q}}}, but this is so that it has the right
form for a composed p.r. function.) 

Suppose now that
\emph{p} codes \uline{1}*\uline{3} so that 
\emph{f\textsubscript{\uline{p}}} expresses the ternary
1\textsuperscript{st}-projection function. It then follows that
\emph{f\textsubscript{\uline{c}}} expresses addition and that \emph{f\textsubscript{\uline{c}}}(\uline{1}, \uline{2})=\uline{3}. The algorithm yields a proof of this by similarly using the information encoded in \emph{c}.

For expediency's sake, the proofs shall utilize inference rules
corresponding to standard evaluation rules for p.r. terms. Where
\uline{\emph{n}}, \uline{{\emph{n}\textsubscript{1}}}, \uline{\emph{n}\textsubscript{2}},~\ldots{} are numerals:

\begin{quote}
(Z) From $\phi(f_{\uline{0}}(\uline{n}))$, $\phi(\uline{0})$ is derivable.

(S) From $\phi(f_{\uline{1}}(\uline{n}))$, $\phi(\uline{n}^\prime)$ is derivable.

(P) If \emph{\uline{c}} codes \emph{\uline{j}}*\emph{\uline{k}} and
$1\!\leq\!\emph{j}\!\leq\!\emph{k}\!\leq\!\emph{c}$, then from
$\phi(f_{\uline{c}}(\uline{n}_{1}, \ldots, \uline{n}_{k}))$, $\phi(\uline{n}_{j})$\\
\phantom{XX.}is derivable.

(C) If \emph{c} codes the string
$\uline{0}\!*\!\uline{\emph{b}}\!*\!\uline{\emph{c}_{1}}\!*\!\ldots*\uline{\emph{c}_{n}}\!*\!\uline{\emph{k}}$, 
and the $\emph{b}\textsuperscript{th}, \emph{c}_{1}\textsuperscript{th}, \ldots$,\\
\phantom{XX.}and
$\emph{c}_{n}\textsuperscript{th}$ axioms define
$f_{\uline{b}}(x_{1}, \ldots, x_{k})$,\!
$f_{\uline{c}_{1}}(x_{1}, \ldots, x_{k})$, \ldots,\\
\phantom{XX.}and
$f_{\uline{c}_{n}}(x_{1}, \ldots, x_{k})$, respectively, then from
$\phi(f_{\uline{c}}(\uline{n}_{1}, \ldots, \uline{n}_{k}))$,\\
\phantom{XX.}$\phi(f_{\uline{b}}(f_{\uline{c}_{1}}(\uline{n}_{1}, \ldots, \uline{n}_{k}), \ldots, f_{\uline{c}_{n}}(\uline{n}_{1}, \ldots, \uline{n}_{k})))$ is derivable.

(R1) If \emph{c} codes the string
`$^\prime$'*\emph{\uline{a}}*\emph{\uline{d}}*\emph{\uline{k}} and the
\emph{a}\textsuperscript{th} and \emph{d}\textsuperscript{th} axioms\\
\phantom{XX.}define $f_{\uline{a}}(x_{1}, \ldots, x_{k})$ and
$f_{\uline{d}}(x_{1}, \ldots, x_{k+2})$, respectively, then\\
\phantom{XX.}from $\phi(f_{\uline{c}}(\uline{n}_{1}, \ldots, \uline{n}_{k}, \uline{0}))$, $\phi(f_{\uline{a}}(\uline{n}_{1}, \ldots, \uline{n}_{k}))$ is derivable.

(R2) Under the same antecedent as (R1), from
$\phi(f_{\uline{c}}(\uline{n}_{1}, \ldots, \uline{n}_{k}, \uline{n}^\prime))$,\\
\phantom{XX.}$\phi(f_{\uline{d}}(f_{\uline{c}}(\uline{n}_{1}, \ldots, \uline{n}_{k}, \uline{n}),
\uline{n}_{1}, \ldots, \uline{n}_{k}, \uline{n}))$ is derivable.

\end{quote}

These are in fact more restrictive than the usual evaluation rules, for
they apply to terms loaded with numerals only. (This is to simplify
things later.) Regardless, we prove in section 5 that the
rules are sound. Soundness means, moreover, that they are mere
shortcuts for what could be proven otherwise in PRA--. Our remarks will therefore bear on PRA--, even though the shortcut rules are not officially part of that system.

We now give an example to illustrate how the subscripts on function symbols guide the application of these rules. To save time, I make use of the standard elimination rule for `\&'; also, some occurrences of `\emph{f\textsubscript{\uline{q}}}( )' are omitted as trivial, although a few are included since they clarify how (R1) is applied. Apart from these omissions, the example is in conformity with the algorithm given later. (To reduce clutter, I omit most underlines, but take heed that numerals appear as such in the proof.) Where $\ulcorner\emph{c}\urcorner\textsuperscript{-1}$ is the string coded by \emph{c}:

\vspace{-6pt}
\begin{footnotesize}
\begin{align*}
1. &\quad f_c(x, 0) = f_q(x)~ \&~ f_c(x, y^{\prime}) = f_d(f_c(x, y), x, y) 
    &&\quad \text{[Axiom for } f\textsubscript{\underline{c}} \text{]} \\
2. &\quad f_c(0^{\prime}, 0) = f_q(0^{\prime})~ \&~ f_c(0^{\prime}, 0^{\prime\prime}) 
    = f_d(f_c(0^{\prime}, 0^{\prime}), 0^{\prime}, 0^{\prime}) 
    &&\quad \text{[(VS), 1: } x/0^{\prime}, y/0^{\prime} \text{]} \\
3. &\quad f_c(0^{\prime}, 0^{\prime\prime}) = f_d(f_c(0^{\prime}, 0^{\prime}), 0^{\prime}, 0^{\prime}) 
    &&\quad \text{[(\&E), 2]} \\
4. &\quad f_c(0^{\prime}, 0^{\prime\prime}) = f_d(f_d(f_c(0^{\prime}, 0), 0^{\prime}, 0), 0^{\prime}, 0^{\prime}) 
        &&\quad \text{[(R2), 3: } \underline{d} = \text{3rd in } \ulcorner c \urcorner^{\text{-1}} \text{]} \\
5. &\quad f_c(0^{\prime}, 0^{\prime\prime}) = f_d(f_d(f_q(0^{\prime}), 0^{\prime}, 0), 0^{\prime}, 0^{\prime}) 
    &&\quad \text{[(R1), 4: } \underline{q}= \text{2nd in } \ulcorner c \urcorner^{\text{-1}} \text{]} \\
6. &\quad f_c(0^{\prime}, 0^{\prime\prime}) = f_d(f_d(0^{\prime}, 0^{\prime}, 0), 0^{\prime}, 0^{\prime}) 
    &&\quad \text{[(P), 5: } q \text{ codes } \underline{1}\!*\!\underline{1} \text{]} \\
7. &\quad f_c(0^{\prime}, 0^{\prime\prime}) = f_d(f_p(f_1(0^{\prime}), 0^{\prime}, 0), 0^{\prime}, 0^{\prime}) 
    &&\quad \text{[(C), 6: } \underline{p} \text{,}~ \underline{1} \text{ are 2nd, 3rd in } \ulcorner d \urcorner^{\text{-}1} \text{]} \\
8. &\quad f_c(0^{\prime}, 0^{\prime\prime}) = f_d(f_p(0^{\prime\prime}, 0^{\prime}, 0), 0^{\prime}, 0^{\prime}) 
    &&\quad \text{[(S), 7: subscript is \uline{1}]} \\
9. &\quad f_c(0^{\prime}, 0^{\prime\prime}) = f_d(0^{\prime\prime}, 0^{\prime}, 0^{\prime}) 
    &&\quad \text{[(P), 8; } p \text{ codes } \underline{1}\!*\!\underline{3} \text{]} \\
10. &\quad f_c(0^{\prime}, 0^{\prime\prime}) = f_p(f_1(0^{\prime\prime}), 0^{\prime}, 0^{\prime}) 
    &&\quad \text{[(C), 9: } \underline{p}\text{,}~ \underline{1} \text{ are 2nd, 3rd in } \ulcorner d \urcorner^{\text{-}1} \text{]} \\
11. &\quad f_c(0^{\prime}, 0^{\prime\prime}) = f_p(0^{\prime\prime\prime}, 0^{\prime}, 0^{\prime}) 
    &&\quad \text{[(S), 10: subscript is \underline{1}]} \\
12. &\quad f_c(0^{\prime}, 0^{\prime\prime}) = 0^{\prime\prime\prime} 
    &&\quad \text{[(P), 11: } p \text{ codes } \underline{1}\!*\!\underline{3}]
\end{align*}

\end{footnotesize}

This example elucidates how subscripts on function symbols enable the construction of the proofs, and in particular, the application of the shortcut rules.

\bigskip
\textbf{4. An algorithm for canonical proofs}

We now present the algorithm for canonical proofs of
\emph{f\textsubscript{\uline{i}}}(\emph{\uline{i}},
\emph{\uline{n}})=\emph{\uline{m}} (when the formula has a proof at
all). As will be detailed later, the algorithm seems free of
unbounded searches, in which case, it would be computable by a p.r.
operation.

Throughout,
the algorithm is often described as operating on linguistic strings rather
than on the codes for these strings. This is for simplicity's sake, and since coding and decoding is p.r., it is of no import.

Terminology: The indices for some functions will code a string
\uline{0}*\emph{\uline{b}}*\uline{\emph{c}\textsubscript{1}}*
\ldots*\emph{\uline{c\textsubscript{n}}}*\emph{\uline{k}}, where the
\emph{b}\textsuperscript{th},
\emph{c}\textsubscript{1}\textsuperscript{th}, \ldots, and
\emph{c\textsubscript{n}}\textsuperscript{th} axioms define
\emph{f\textsubscript{\uline{b}}}(\emph{x}\textsubscript{1}, \ldots,
\emph{x\textsubscript{k}}),
\emph{f\textsubscript{\uline{c}}}\textsubscript{\uline{1}}(\emph{x}\textsubscript{1},
\ldots, \emph{x\textsubscript{k}}), \ldots, and
\emph{f\textsubscript{\uline{cn}}}(\emph{x}\textsubscript{1}, \ldots,
\emph{x\textsubscript{k}}), respectively. Call such an index a
``composition index.'' Other indices will code a string
`$^\prime$'*\emph{\uline{a}}*\emph{\uline{d}}, where the
\emph{a}\textsuperscript{th} and \emph{d}\textsuperscript{th} axioms
define \emph{f\textsubscript{\uline{a}}}(\emph{x}\textsubscript{1},
\ldots, \emph{x\textsubscript{k}}) and
\emph{f\textsubscript{\uline{d}}}(\emph{x}\textsubscript{1}, \ldots,
\emph{x\textsubscript{k}}\textsubscript{+2}), respectively. Call such an
index a ``recursion index.'' Also, let us jointly refer to composition
indices and recursion indices as ``complex indices.'' The other
indices we shall call ``simple indices.''

Suppose now that \emph{e}(\emph{i})
enumerates the axiom for \(f_{\underline{i}}^{\underline{2}}\). Then, \emph{s}(\emph{i},\emph{n}) is the code of the sequence that is determined as follows:

\bigskip
\textbf{Step 1.} {[}Basic Binary Projections{]} If \emph{i} does not code \emph{\uline{j}}*\uline{2}, where \emph{j} equals 1 or 2, then go to Step 2. Otherwise:
  
\begin{enumerate}[label=\Alph*.]
  \item
If \emph{j}=1, output the following sequence and then halt:
  \end{enumerate}

\begin{longtable}[]{@{}
  >{\raggedright\arraybackslash}p{(\columnwidth - 2\tabcolsep) * \real{0.7}}
  >{\raggedright\arraybackslash}p{(\columnwidth - 2\tabcolsep) * \real{0.3}}@{}}
\hspace{.5cm}1. \emph{f\textsubscript{\uline{i}}}(\emph{x}, \emph{y})=\emph{x} &
{[}Axiom for \emph{f\textsubscript{\uline{i}}}{]} \\
\hspace{.5cm}2. \emph{f\textsubscript{\uline{i}}}(\emph{\uline{i}},
\emph{\uline{n}})=\emph{\uline{i}} & {[}(VS), 1:
\emph{x}/\emph{\uline{i}}, \emph{y}/\emph{\uline{n}}{]} \\

\end{longtable}

\begin{enumerate}
\def\labelenumi{\Alph{enumi}.}
\setcounter{enumi}{1}
\item
  If \emph{j}=2, output the following sequence and then halt:
\end{enumerate}

\begin{longtable}[]{@{}
  >{\raggedright\arraybackslash}p{(\columnwidth - 2\tabcolsep) * \real{0.7}}
  >{\raggedright\arraybackslash}p{(\columnwidth - 2\tabcolsep) * \real{0.3}}@{}}

\hspace{.5cm}1. \emph{f\textsubscript{\uline{i}}}(\emph{x}, \emph{y})=\emph{y} &
{[}Axiom for \emph{f\textsubscript{\uline{i}}}{]} \\
\hspace{.5cm}2. \emph{f\textsubscript{\uline{i}}}(\emph{\uline{i}},
\emph{\uline{n}})=\emph{\uline{n}} & {[}(VS), 1:
\emph{x}/\emph{\uline{i}}, \emph{y}/\emph{\uline{n}}{]} \\
\end{longtable}

\textbf{Step 2.} {[}Getting the Baseline Term for P.R. Recursions{]} If \emph{i} is not a recursion index, then go to Step 3. Otherwise:

  \begin{enumerate}[label=\Alph*.]
  \def\labelenumii{\Alph{enumii}.}
  \item
    If \emph{n}=0, start the sequence as follows and then go to Step 4:
  \end{enumerate}

\begin{longtable}[]{@{}
  >{\raggedright\arraybackslash}p{(\columnwidth - 2\tabcolsep) * \real{0.7}}
  >{\raggedright\arraybackslash}p{(\columnwidth - 2\tabcolsep) * \real{0.3}}@{}}

\hspace{.5cm}1. \emph{f\textsubscript{\uline{i}}}(\emph{x},
0)=\emph{f\textsubscript{\uline{a}}}(\emph{x}) \&
\emph{f\textsubscript{\uline{i}}}(\emph{x},
\emph{y}$^\prime$)=\emph{f\textsubscript{\uline{d}}}(\emph{f\textsubscript{\uline{i}}}(\emph{x},
\emph{y}), \emph{x}, \emph{y}) & {[}Axiom for
\emph{f\textsubscript{\uline{i}}}{]} \\
\hspace{.5cm}2. \emph{f\textsubscript{\uline{i}}}(\emph{\uline{i}},
0)=\emph{f\textsubscript{\uline{a}}}(\emph{\uline{i}}) \&
\emph{f\textsubscript{\uline{i}}}(\emph{\uline{i}},
\emph{\uline{n}}$^\prime$)=\emph{f\textsubscript{\uline{d}}}(\emph{f\textsubscript{\uline{i}}}(\emph{\uline{i}},
\emph{\uline{n}}), \emph{\uline{i}}, \emph{\uline{n}}) & {[}(VS), 1:
\emph{x}/\emph{\uline{i}}, \emph{y}/\emph{\uline{n}}{]} \\
\hspace{.5cm}3. \emph{f\textsubscript{\uline{i}}}(\emph{\uline{i}},
0)=\emph{f\textsubscript{\uline{a}}}(\emph{\uline{i}}) & {[}(\&E),
2{]}\\
\end{longtable}

\begin{enumerate}
\def\labelenumi{\Alph{enumi}.}
\setcounter{enumi}{1}
\item
  If $\emph{n}\neq0$, start the sequence as follows and then go to Step 4:
\end{enumerate}

\begin{longtable}[]{@{}
  >{\raggedright\arraybackslash}p{(\columnwidth - 2\tabcolsep) * \real{0.7}}
  >{\raggedright\arraybackslash}p{(\columnwidth - 2\tabcolsep) * \real{0.3}}@{}}

\hspace{.5cm}1. \emph{f\textsubscript{\uline{i}}}(\emph{x},
0)=\emph{f\textsubscript{\uline{a}}}(\emph{x}) \&
\emph{f\textsubscript{\uline{i}}}(\emph{x},
\emph{y}$^\prime$)=\emph{f\textsubscript{\uline{d}}}(\emph{f\textsubscript{\uline{i}}}(\emph{x},
\emph{y}), \emph{x}, \emph{y}) & {[}Axiom for
\emph{f\textsubscript{\uline{i}}}{]} \\
\hspace{.5cm}2. \emph{f\textsubscript{\uline{i}}}(\emph{\uline{i}},
0)=\emph{f\textsubscript{\uline{a}}}(\emph{\uline{i}}) \&
\emph{f\textsubscript{\uline{i}}}(\emph{\uline{i}},
\uline{\emph{n}-1}$^\prime$)=\emph{f\textsubscript{\uline{d}}}(\emph{f\textsubscript{\uline{i}}}(\emph{\uline{i}},
\uline{\emph{n}-1}), \emph{\uline{i}}, \uline{\emph{n}-1}) & {[}(VS), 1:
\emph{x}/\emph{\uline{i}}, \emph{y}/\uline{\emph{n}-1}{]} \\
\hspace{.5cm}3. \emph{f\textsubscript{\uline{i}}}(\emph{\uline{i}},
\uline{\emph{n}-1}$^\prime$)=\emph{f\textsubscript{\uline{d}}}(\emph{f\textsubscript{\uline{i}}}(\emph{\uline{i}},
\uline{\emph{n}-1}), \emph{\uline{i}}, \uline{\emph{n}-1}) & {[}(\&E),
2{]}\\
\end{longtable}

\textbf{Step 3.} {[}Getting the Baseline Term for Other Functions{]} Start the sequence as follows and then go to Step 4:

\begin{longtable}[]{@{}
  >{\raggedright\arraybackslash}p{(\columnwidth - 2\tabcolsep) * \real{0.7}}
  >{\raggedright\arraybackslash}p{(\columnwidth - 2\tabcolsep) * \real{0.3}}@{}}

\hspace{.5cm}1. \emph{f\textsubscript{\uline{i}}}(\emph{x},
\emph{y})=$\phi$\emph{\textsubscript{i}}(\emph{x}, \emph{y}) & {[}Axiom for
\emph{f\textsubscript{\uline{i}}}{]} \\
\hspace{.5cm}2. \emph{f\textsubscript{\uline{i}}}(\emph{\uline{i}},
\emph{\uline{n}})=$\phi$\emph{\textsubscript{i}}(\emph{\uline{i},
\uline{n}}) & {[}(VS), 1: \emph{x}/\emph{\uline{i}},
\emph{y}/\emph{\uline{n}}{]} \\
\end{longtable}

\textbf{Step 4.} {[}Initializing the Counter{]} The latest line of the sequence is an equality; call the right-hand term the ``baseline term'' for the sequence. Count how many times `\emph{f}' occurs in the baseline term and initialize the counter at that number (not greater than the total number of symbols in the wff\footnote{Cf. Gödel's function \emph{l}(\emph{x}); \#7 in his list of p.r. functions and relations, p. 182.}). Go to Step 5.
\\

 \textbf{Step 5.} {[}Entering the Main Loop.{]} If the counter is at zero, halt. Otherwise, (re)start the Main Loop: For the present iteration of the Loop, let $\epsilon$ be the equation on the latest line. To the right of `=', find the nn-term embedded in the most parentheses (a.k.a., the ``innermost'' nn-term). When there is a tie, choose the rightmost one. Let $\tau$ be the rightmost, innermost nn-term (to the right of `=') for the current iteration of the Loop. (N.B., $\tau$ will be loaded with numerals only; we prove this later.) Go to Step 6.
\\

\textbf{Step 6.} {[}Computing $\tau${]} Check the index \emph{c} for $\tau$. It is either 0, 1,
  or codes \emph{\uline{j}}*\emph{\uline{k}}, where
  $1\!\leq\!\emph{j}\!\leq\!\emph{k}\!\leq\!\emph{c}$---alternatively, it is a composition or
  recursion index.
\vspace{6pt}

\begin{enumerate}
\def\labelenumi{\Alph{enumi}.}
\item
  If \emph{c}=0, apply (Z): Add a line where $\tau$ is replaced in $\epsilon$ with
  \uline{0}. Subtract 1 from the counter and go back to Step 5.
\item
  If \emph{c}=1, then $\tau$ is loaded with some \emph{\uline{m}}. Apply (S):
  Add a line where $\tau$ is replaced in $\epsilon$ with \emph{\uline{m}}$^\prime$. Subtract 1
  from the counter and go back to Step 5.
\item
  If \emph{c} codes \emph{\uline{j}}*\emph{\uline{k}} where
  $1\!\leq\!\emph{j}\!\leq\!\emph{k}\!\leq\!\emph{c}$, then apply (P): Add a line where $\tau$ is
  replaced in $\epsilon$ with the numeral in the \emph{j}\textsuperscript{th}
  position of $\tau$. Subtract 1 from the counter and go back to Step 5.
\item
  If \emph{c} is a composition index, then apply (C): Add a line where $\tau$
  is replaced in $\epsilon$ with
  \emph{f\textsubscript{\uline{b}}}(\emph{f\textsubscript{\uline{c}}}\textsubscript{\uline{1}}(\uline{\emph{n}\textsubscript{1}},
  \ldots, \emph{\uline{n\textsubscript{k}}}), \ldots,
  \emph{f\textsubscript{\uline{cl}}}(\uline{\emph{n}\textsubscript{1}},
  \ldots, \emph{\uline{n\textsubscript{k}}})), where
  \uline{\emph{n}\textsubscript{1}}, \ldots,
  \emph{\uline{n\textsubscript{k}}} are the same as in $\tau$, and
  \emph{\uline{b}} is the 2\textsuperscript{nd} member of the string
  coded by \emph{c}, \uline{\emph{c}\textsubscript{1}} is
  3\textsuperscript{rd} member of the coded string coded by \emph{c},
  \ldots, and \emph{\uline{c\textsubscript{l}}} is
  \emph{l}+2\textsuperscript{th} member of the string coded by \emph{c}.
  Add \emph{l} to the counter and go back to Step 5.

\item
  If \emph{c} is a recursion index, then:

  \begin{enumerate}
  \def\labelenumii{\roman{enumii}.}
  \item
    
    If $\tau$ is loaded with \uline{\emph{n}\textsubscript{1}}, \ldots,
    \emph{\uline{n\textsubscript{k}}}, \uline{0}, then apply (R1): Add a
    line where $\tau$ is replaced in $\epsilon$ with
    \emph{f\textsubscript{\uline{a}}}(\uline{\emph{n}\textsubscript{1}},
    \ldots, \emph{\uline{n\textsubscript{k}}}), where \emph{\uline{a}}
    is the 2\textsuperscript{nd} member of the string coded by \emph{c}.
    Leave the counter unchanged and go back to Step 5.
  
  \item
    
    If $\tau$ is loaded with \uline{\emph{n}\textsubscript{1}}, \ldots,
    \emph{\uline{n\textsubscript{k}}}, \emph{\uline{m}}$\neq$\uline{0}, then apply (R2): Add a line where $\tau$ is
    replaced in $\epsilon$ with
    \emph{f\textsubscript{\uline{d}}}(\emph{f\textsubscript{\uline{c}}}(\uline{\emph{n}\textsubscript{1}},
    \ldots, \emph{\uline{n\textsubscript{k}}}, \uline{\emph{m}-1}),
    \uline{\emph{n}\textsubscript{1}}, \ldots,
    \emph{\uline{n\textsubscript{k}}}, \uline{\emph{m}--1}), where
    \emph{\uline{c}} is the same as in $\tau$ and \emph{\uline{d}} is the
    3\textsuperscript{rd} member of the string coded by \emph{c}. Add 1
    to the counter and go back to Step 5.
   
  \end{enumerate}
\end{enumerate}

Again, the algorithm generates a proof where
\emph{f\textsubscript{\uline{i}}}(\emph{\uline{i}}, \emph{\uline{n}}) is
computed for any \emph{n}, provided that the
\emph{i}\textsuperscript{th} axiom for the function symbols defines a
binary symbol. The next section proves this, and in the
section after that, an argument is given for why the algorithm is
p.r.

\bigskip
\textbf{5. Proof that the algorithm is correct}

In precise terms, the claim to be proved is:

\begin{quote}

\textbf{Theorem}: If \emph{e}(\emph{i}) codes the axiom defining
\(f_{\underline{i}}^{\underline{2}}\), then given \emph{i} and any \emph{n}, there is an \emph{m}
such that the algorithm generates a unique proof of\\
\(f_{\underline{i}}\)(\emph{\uline{i}},
\emph{\uline{n}})=\emph{\uline{m}}, at which point the algorithm halts.

\end{quote}

Note that here and elsewhere in this section, our statements employ
unbounded quantification; however, this alone does not undermine that
the algorithm is p.r. \textbf{Theorem} is a claim about the algorithm, not
part of the algorithm itself.

Establishing \textbf{Theorem} is best approached by considering p.r. terms
with simple indices first, and then considering separately those with
complex indices.

Simple Indices: Vacuously, \textbf{Theorem} holds if the index \emph{i}=0;
after all, the 0\textsuperscript{th} function symbol is not binary, and
a binary function-symbol with subscript `0' has no axiom in the
system. Similar remarks apply with every other simple index, except an index that codes  \uline{1}* \uline{2} or \uline{2}* \uline{2}. But in those cases, the relevant function is a basic binary projection, and the
algorithm clearly handles these under Step 1 in a way that satisfies
\textbf{Theorem}. 

Complex Indices: We want to show that \textbf{Theorem} holds when
\emph{i} is a complex index. This part of the proof will require four
lemmas. The first is as follows:

\begin{quote}

\uline{Lemma 1}: If \(f_{\underline{i}}^{\underline{2}}\) has a complex index,
then for any \emph{n}, the algorithm starts with a \emph{proof} of
\emph{f\textsubscript{\uline{i}}}(\emph{\uline{i}}, \emph{\uline{n}})=$\beta$,
where $\beta$ is a baseline term and has no free variables.

\end{quote}

Proof: If \emph{i} is complex, there are three kinds of case to
consider. The first two are where \emph{i} is a recursion index and where \emph{n}=0 and $\emph{n}\neq0$, respectively. The third
is where \emph{i} is a composition index. In the first case, the
algorithm starts:

\begin{longtable}[]{@{}
  >{\raggedright\arraybackslash}p{(\columnwidth - 2\tabcolsep) * \real{0.7}}
  >{\raggedright\arraybackslash}p{(\columnwidth - 2\tabcolsep) * \real{0.3}}@{}}
\hspace{.5cm}1. \emph{f\textsubscript{\uline{i}}}(\emph{x},
0)=\emph{f\textsubscript{\uline{a}}}(\emph{x}) \&
\emph{f\textsubscript{\uline{i}}}(\emph{x},
\emph{y}$^\prime$)=\emph{f\textsubscript{\uline{b}}}(\emph{f\textsubscript{\uline{i}}}(\emph{x},
\emph{y}), \emph{x}, \emph{y}) & {[}Axiom for
\emph{f\textsubscript{\uline{i}}}{]} \\
\hspace{.5cm}2. \emph{f\textsubscript{\uline{i}}}(\emph{\uline{i}},
0)=\emph{f\textsubscript{\uline{a}}}(\emph{\uline{i}}) \&
\emph{f\textsubscript{\uline{i}}}(\emph{\uline{i}},
\emph{\uline{n}}$^\prime$)=\emph{f\textsubscript{\uline{b}}}(\emph{f\textsubscript{\uline{i}}}(\emph{\uline{i}},
\emph{\uline{n}}), \emph{\uline{i}}, \emph{\uline{n}}) & {[}(VS), 1:
\emph{x}/\emph{\uline{i}}, \emph{y}/\emph{\uline{n}}{]} \\
\hspace{.5cm}3. \emph{f\textsubscript{\uline{i}}}(\emph{\uline{i}},
0)=\emph{f\textsubscript{\uline{a}}}(\emph{\uline{i}}) & {[}(\&E),
2{]} \\
\end{longtable}
\vspace{-6pt}
The rules used are obviously sound; also,
\emph{f\textsubscript{\uline{a}}}(\emph{\uline{i}}) has no free
variables and is by definition a baseline term. So the third line
verifies \uline{Lemma 1} in this first case. In the second case, the
algorithm begins:

\begin{longtable}[]{@{}
  >{\raggedright\arraybackslash}p{(\columnwidth - 2\tabcolsep) * \real{0.7}}
  >{\raggedright\arraybackslash}p{(\columnwidth - 2\tabcolsep) * \real{0.3}}@{}}

\hspace{.5cm}1. \emph{f\textsubscript{\uline{i}}}(\emph{x},
0)=\emph{f\textsubscript{\uline{a}}}(\emph{x}) \&
\emph{f\textsubscript{\uline{i}}}(\emph{x},
\emph{y}$^\prime$)=\emph{f\textsubscript{\uline{b}}}(\emph{f\textsubscript{\uline{i}}}(\emph{x},
\emph{y}), \emph{x}, \emph{y}) & {[}Axiom for
\emph{f\textsubscript{\uline{i}}}{]} \\
\hspace{.5cm}2. \emph{f\textsubscript{\uline{i}}}(\emph{\uline{i}},
0)=\emph{f\textsubscript{\uline{a}}}(\emph{\uline{i}}) \&
\emph{f\textsubscript{\uline{i}}}(\emph{\uline{i}},
\uline{\emph{n}-1}$^\prime$)=\emph{f\textsubscript{\uline{b}}}(\emph{f\textsubscript{\uline{i}}}(\emph{\uline{i}},
\uline{\emph{n}-1}), \emph{\uline{i}}, \uline{\emph{n}-1}) & {[}(VS), 1:
\emph{x}/\emph{\uline{i}}, \emph{y}/\uline{\emph{n}-1}{]} \\
\hspace{.5cm}3. \emph{f\textsubscript{\uline{i}}}(\emph{\uline{i}},
\uline{\emph{n}-1}$^\prime$)=\emph{f\textsubscript{\uline{b}}}(\emph{f\textsubscript{\uline{i}}}(\emph{\uline{i}},
\uline{\emph{n}-1}), \emph{\uline{i}}, \uline{\emph{n}-1}) & {[}(\&E),
2{]} \\
\end{longtable}
\vspace{-6pt}
The same rules are used as before, and again,
\emph{f\textsubscript{\uline{b}}}(\emph{f\textsubscript{\uline{i}}}(\emph{\uline{i}},
\uline{\emph{n}-1}), \emph{\uline{i}}, \uline{\emph{n}-1}) has no free
variables and is by definition a baseline term. So the third line
verifies \uline{Lemma 1} in the second case. In the third case, the
algorithm starts as follows:

\begin{longtable}[]{@{}
  >{\raggedright\arraybackslash}p{(\columnwidth - 2\tabcolsep) * \real{0.7}}
  >{\raggedright\arraybackslash}p{(\columnwidth - 2\tabcolsep) * \real{0.3}}@{}}

\hspace{.5cm}1. \emph{f\textsubscript{\uline{i}}}(\emph{x},
\emph{y})=$\phi$\emph{\textsubscript{i}}(\emph{x}, \emph{y}) & {[}Axiom for
\emph{f\textsubscript{\uline{i}}}{]} \\
\hspace{.5cm}2. \emph{f\textsubscript{\uline{i}}}(\emph{\uline{i}},
\emph{\uline{n}})=$\phi$\emph{\textsubscript{i}}(\emph{\uline{i},
\uline{n}}) & {[}(VS), 1: \emph{x}/\emph{\uline{i}},
\emph{y}/\emph{\uline{n}}{]} \\
\end{longtable}

\vspace{-6pt}
The rule is sound; $\phi$\emph{\textsubscript{i}}(\emph{\uline{i},
\uline{n}}) has no free variables and is by definition a baseline term.
So the second line verifies \uline{Lemma 1} in the third case, which
completes the proof of \uline{Lemma 1}.

\emph{Remark}: Like the previous two cases, the baseline term in the
third case will have at least one occurrence of `\emph{f}'. After all,
\emph{f\textsubscript{\uline{i}}} in the third case expresses a composed
p.r. function, and the axiom for \emph{f\textsubscript{\uline{i}}}
therefore will contain at least one `\emph{f}' on the right-hand side
(and the baseline term is an instantiation of the right-hand side).

The second lemma required is the following:
\begin{quote}
\uline{Lemma 2}: {[}Soundness{]} If the algorithm applies (Z), (S), (P),
(C), (R1), or (R2) to a line in order to produce a new line of a
sequence, then in the standard model, if the former line is true, so is
the latter.
\end{quote}

The lemma is claiming that the algorithm uses the shortcut rules in a
way that is sound in the standard model. We prove this by considering
its use of each of the rules.

Preliminary observation: The algorithm is designed to apply the shortcut
rules only to the rightmost, innermost nn-term $\tau$ on a line. So, $\tau$ will
be loaded with numerals only. For if $\tau$ were loaded with a nn-term, it
would not be the innermost. (Also, we also know from \uline{Lemma 1}
that a baseline term has no free variables, and none of the inference
rules introduce free variables into the picture.)

\emph{Applying (Z)}: Suppose the algorithm applies (Z) so to produce
\emph{f\textsubscript{\uline{i}}}(\emph{\uline{i}},
\emph{\uline{n}})=$\phi$($\tau$\textsubscript{2}) from
\emph{f\textsubscript{\uline{i}}}(\emph{\uline{i}},
\emph{\uline{n}})=$\phi$($\tau$\textsubscript{1}). Then, per the instructions on
Step 6A, $\tau$\textsubscript{1} must be of the form
\emph{f}\textsubscript{\uline{0}}(\emph{\uline{m}}), for some numeral
\emph{\uline{m}}. Also, per those instructions, $\tau$\textsubscript{2} must
be \uline{0}. Since \emph{f}\textsubscript{\uline{0}} expresses the
constantly-zero function, \uline{Lemma 2} is verified in the case of
(Z).

\emph{Applying (S)}: Suppose the algorithm applies (S) so to produce
\emph{f\textsubscript{\uline{i}}}(\emph{\uline{i}},
\emph{\uline{n}})=$\phi$($\tau$\textsubscript{2}) from
\emph{f\textsubscript{\uline{i}}}(\emph{\uline{i}},
\emph{\uline{n}})=$\phi$($\tau$\textsubscript{1}). Then, per the instructions on
Step 6B, $\tau$\textsubscript{1} must be of the form
\emph{f}\textsubscript{\uline{1}}(\emph{\uline{m}}), for some
\emph{\uline{m}}. Also, per those instructions, $\tau$\textsubscript{2} must
be \emph{\uline{m}}$^\prime$. Thus, since \emph{f}\textsubscript{\uline{1}}
expresses the successor function, \uline{Lemma 2} is verified in the
case of (S).

\emph{Applying (P)}: Suppose the algorithm applies (P) so to produce
\emph{f\textsubscript{\uline{i}}}(\emph{\uline{i}},
\emph{\uline{n}})=$\phi$($\tau$\textsubscript{2}) from
\emph{f\textsubscript{\uline{i}}}(\emph{\uline{i}},
\emph{\uline{n}})=$\phi$($\tau$\textsubscript{1}) Then, per the instructions on
Step 6C, $\tau$\textsubscript{1} must be of the form
\emph{f\textsubscript{\uline{c}}}(\uline{\emph{n}\textsubscript{1}},
\ldots, \emph{\uline{n\textsubscript{k}}}), where \emph{c} codes
\emph{\uline{j}}*\emph{\uline{k}} and $1\!\leq\!\emph{j}\!\leq\!\emph{k}\!\leq\!\emph{c}$.
Also, per those instructions, $\tau$\textsubscript{2} must be
\emph{\uline{n\textsubscript{j}}}. Thus, since
\emph{f\textsubscript{\uline{c}}} expresses the \emph{k}-ary
\emph{j}\textsuperscript{th}-projection function, \uline{Lemma 2} is
verified in the case of (P).

\emph{Applying (C)}: Suppose the algorithm applies (C) so to produce
\emph{f\textsubscript{\uline{i}}}(\emph{\uline{i}},
\emph{\uline{n}})=$\phi$($\tau$\textsubscript{2}) from
\emph{f\textsubscript{\uline{i}}}(\emph{\uline{i}},
\emph{\uline{n}})=$\phi$($\tau$\textsubscript{1}). Then, per the instructions on
Step 6D, $\tau$\textsubscript{1} must be of the form
\emph{f\textsubscript{\uline{c}}}(\uline{\emph{n}\textsubscript{1}},
\ldots, \emph{\uline{n\textsubscript{k}}}), where \emph{c} codes
\uline{0}*\emph{\uline{b}}*\uline{\emph{c}\textsubscript{1}}*\ldots*\emph{\uline{c\textsubscript{l}}}*\emph{\uline{k}.}
Also, per those instructions, $\tau$\textsubscript{2} must be
\emph{f\textsubscript{\uline{b}}}(\emph{f\textsubscript{\uline{c}}}\textsubscript{\uline{1}}(\uline{\emph{n}\textsubscript{1}},
\ldots, \emph{\uline{n\textsubscript{k}}}), \ldots,
\emph{f\textsubscript{\uline{cl}}}(\uline{\emph{n}\textsubscript{1}},
\ldots, \emph{\uline{n\textsubscript{k}}})). Observe that, given how the
indices for composed functions are assigned,
\emph{f\textsubscript{\uline{c}}} expresses a function where the
\emph{l}-ary function expressed by \emph{f\textsubscript{\uline{b}}} is
composed with the \emph{k}-ary functions expressed by
\emph{f\textsubscript{\uline{c}}}\textsubscript{\uline{1}}, \ldots, and
\emph{f\textsubscript{\uline{cl}}}. The consequence is that
\emph{f\textsubscript{\uline{c}}}(\uline{\emph{n}\textsubscript{1}},
\ldots, \emph{\uline{n\textsubscript{k}}}) is co-referential with\\
\emph{f\textsubscript{\uline{b}}}(\emph{f\textsubscript{\uline{c}}}\textsubscript{\uline{1}}(\uline{\emph{n}\textsubscript{1}},
\ldots, \emph{\uline{n\textsubscript{k}}}), \ldots,
\emph{f\textsubscript{\uline{cl}}}(\uline{\emph{n}\textsubscript{1}},
\ldots, \emph{\uline{n\textsubscript{k}}})). Thus, the transition from
\emph{f\textsubscript{\uline{i}}}(\emph{\uline{i}},
\emph{\uline{n}})=$\phi$($\tau$\textsubscript{1}) to
\emph{f\textsubscript{\uline{i}}}(\emph{\uline{i}},
\emph{\uline{n}})=$\phi$($\tau$\textsubscript{2}) is sound, and \uline{Lemma 2} is
verified in the case of (C).

\emph{Applying (R1)}: Suppose the algorithm applies (R1) so to produce
\emph{f\textsubscript{\uline{i}}}(\emph{\uline{i}},
\emph{\uline{n}})=$\phi$($\tau$\textsubscript{2}) from
\emph{f\textsubscript{\uline{i}}}(\emph{\uline{i}},
\emph{\uline{n}})=$\phi$($\tau$\textsubscript{1}). Then, per the instructions on
Step 6Ei, $\tau$\textsubscript{1} must be of the form
\emph{f\textsubscript{\uline{c}}}(\uline{\emph{n}\textsubscript{1}},
\ldots, \emph{\uline{n\textsubscript{k}}}, \uline{0}) where \emph{c}
codes `$^\prime$'*\emph{\uline{a}}*\emph{\uline{d}}*\emph{\uline{k}.} Also, per
those instructions, $\tau$\textsubscript{2} must be
\emph{f\textsubscript{\uline{a}}}(\uline{\emph{n}\textsubscript{1}},
\ldots, \emph{\uline{n\textsubscript{k}}}). Observe that, given how the
indices for composed functions are assigned,
\emph{f\textsubscript{\uline{c}}} expresses a function which is
recursively defined by \emph{f\textsubscript{\uline{a}}} and
\emph{f\textsubscript{\uline{d}}}. The consequence is that
\emph{f\textsubscript{\uline{c}}}(\uline{\emph{n}\textsubscript{1}},
\ldots, \emph{\uline{n\textsubscript{k}}}, \uline{0}) co-refers with
\emph{f\textsubscript{\uline{a}}}(\uline{\emph{n}\textsubscript{1}},
\ldots, \emph{\uline{n\textsubscript{k}}}). So the transition from
\emph{f\textsubscript{\uline{i}}}(\emph{\uline{i}},
\emph{\uline{n}})=$\phi$($\tau$\textsubscript{1}) to
\emph{f\textsubscript{\uline{i}}}(\emph{\uline{i}},
\emph{\uline{n}})=$\phi$($\tau$\textsubscript{2}) is sound, and \uline{Lemma 2} is
verified in the case of (R1).

\emph{Applying (R2)}: Suppose the algorithm applies (R2) so to produce
\emph{f\textsubscript{\uline{i}}}(\emph{\uline{i}},
\emph{\uline{n}})=$\phi$($\tau$\textsubscript{2}) from
\emph{f\textsubscript{\uline{i}}}(\emph{\uline{i}},
\emph{\uline{n}})=$\phi$($\tau$\textsubscript{1}). Then, per the instructions on
Step 6Eii, $\tau$\textsubscript{1} must be of the form
\emph{f\textsubscript{\uline{c}}}(\uline{\emph{n}\textsubscript{1}},
\ldots, \emph{\uline{n\textsubscript{k}}}, \emph{\uline{m}}) for
\emph{m}\textgreater0, where \emph{c} codes
`$^\prime$'*\emph{\uline{a}}*\emph{\uline{d}}*\emph{\uline{k}.} Also, per those
instructions, $\tau$\textsubscript{2} must be
\emph{f\textsubscript{\uline{d}}}(\emph{f\textsubscript{\uline{c}}}(\uline{\emph{n}\textsubscript{1}},
\ldots, \emph{\uline{n\textsubscript{k}}}, \uline{\emph{m}-1}),
\uline{\emph{n}\textsubscript{1}}, \ldots,
\emph{\uline{n\textsubscript{k}}}, \uline{\emph{m}--1}). Observe that,
given how the indices for composed functions are assigned,
\emph{f\textsubscript{\uline{c}}} expresses a function which is
recursively defined by \emph{f\textsubscript{\uline{a}}} and
\emph{f\textsubscript{\uline{d}}}. The consequence is that
\emph{f\textsubscript{\uline{c}}}(\uline{\emph{n}\textsubscript{1}},
\ldots, \emph{\uline{n\textsubscript{k}}}, \emph{\uline{m}}) co-refers
with
\emph{f\textsubscript{\uline{d}}}(\emph{f\textsubscript{\uline{c}}}(\uline{\emph{n}\textsubscript{1}},
\ldots, \emph{\uline{n\textsubscript{k}}}, \uline{\emph{m}-1}),
\uline{\emph{n}\textsubscript{1}}, \ldots,
\emph{\uline{n\textsubscript{k}}}, \uline{\emph{m}--1}). So the
transition from \emph{f\textsubscript{\uline{i}}}(\emph{\uline{i}},
\emph{\uline{n}})=$\phi$($\tau$\textsubscript{1}) to
\emph{f\textsubscript{\uline{i}}}(\emph{\uline{i}},
\emph{\uline{n}})=$\phi$($\tau$\textsubscript{2}) is sound, and \uline{Lemma 2} is
verified in the case of (R2). This completes the proof of \uline{Lemma
2}.

The third lemma will be needed for proving the fourth lemma, later:

\begin{quote}
\uline{Lemma 3}: If $\beta$ is the baseline term for a sequence and the
algorithm produces a derivation from
\emph{f\textsubscript{\uline{i}}}(\emph{\uline{i}}, \emph{\uline{n}})=$\beta$
to \emph{f\textsubscript{\uline{i}}}(\emph{\uline{i}},
\emph{\uline{n}})=$\gamma$, then at that point, the counter is set to the
number of times `\emph{f}' occurs in $\gamma$.
\end{quote}

We argue this by induction on the number of times that the Main Loop
iterates in deriving \emph{f\textsubscript{\uline{i}}}(\emph{\uline{i}},
\emph{\uline{n}})=$\gamma$ from
\emph{f\textsubscript{\uline{i}}}(\emph{\uline{i}}, \emph{\uline{n}})=$\beta$.

\emph{Base Case}: The Main Loop iterates zero times. Then, the
derivation just consists in the single line
\emph{f\textsubscript{\uline{i}}}(\emph{\uline{i}}, \emph{\uline{n}})=$\beta$,
meaning that $\gamma$=$\beta$. Moreover, the counter is initially set to the number
of `\emph{f}'s in $\beta$ by stipulation. So the base case verifies
\uline{Lemma 3}.

\emph{Inductive Cases}: Suppose that if the algorithm iterates
the Main Loop \emph{m} times to produce a derivation from
\emph{f\textsubscript{\uline{i}}}(\emph{\uline{i}}, \emph{\uline{n}})=$\beta$,
from \emph{f\textsubscript{\uline{i}}}(\emph{\uline{i}},
\emph{\uline{n}})=$\delta$, then and at that point, the counter matches the
number of `\emph{f}'s in $\delta$. We want to show that if the algorithm
iterates the Main Loop one more time to derive
\emph{f\textsubscript{\uline{i}}}(\emph{\uline{i}}, \emph{\uline{n}})=$\gamma$,
the counter then matches the number of `\emph{f}'s in $\gamma$. There are six
cases to consider.

\emph{Case 1}: The \emph{m}+1\textsuperscript{th} iteration of the Loop applies (Z). Then, the algorithm adds a line
which is the same as the previous one, except that the rightmost, innermost nn-term
in $\delta$ is replaced with \uline{0}. The replacement means that the right
side of the newest line has one less occurrence of `\emph{f}'.
Similarly, whereas the counter was set to the number of `\emph{f}'s in
$\delta$, the algorithm now reduces it by 1. Thus, the counter matches the
number of `\emph{f}'s to the right of `=' in the newest line. So,
\emph{Case 1} verifies \uline{Lemma 3}.

\emph{Case 2}: The \emph{m}+1\textsuperscript{th} iteration of the Loop applies (S). Then, the algorithm adds a line
which is the same as the previous one, except that the rightmost, innermost nn-term
in $\delta$, which is of the form \emph{f}\textsubscript{1}(\emph{\uline{m}}),
is replaced with \emph{\uline{m}}$^\prime$. The replacement means that the
newest line has one less occurrence of `\emph{f}'. Similarly, whereas
the counter was set to the number of `\emph{f}'s in $\delta$, the algorithm now
reduces it by 1. Thus, the counter matches the number of `\emph{f}'s to
the right of `=' in the newest line. So, \emph{Case 2} verifies
\uline{Lemma 3}.

\emph{Case 3}: The \emph{m}+1\textsuperscript{th} iteration of the Loop applies (P). Then, the algorithm adds a line
which is the same as the previous one, except that the rightmost, innermost nn-term
in $\delta$, which is of the form
\emph{f\textsubscript{\uline{c}}}(\emph{\uline{n\textsubscript{1}}},
\ldots, \emph{\uline{n\textsubscript{k}}}), has been replaced by a
numeral \emph{\uline{n\textsubscript{j}}}, where $\emph{j}\!\leq\!\emph{k}$. The
replacement means that the right side of the newest line has one less
occurrence of `\emph{f}'. Similarly, whereas the counter was set to the
number of `\emph{f}'s in $\delta$, the algorithm now reduces it by 1. Thus, the
counter matches the number of `\emph{f}'s to the right of `=' in the
newest line. So, \emph{Case 3} verifies \uline{Lemma 3}.

\emph{Case 4}: The \emph{m}+1\textsuperscript{th} iteration of the Loop applies (C). In such a case, the algorithm adds a line
which is the same as the previous one, except that the rightmost, innermost nn-term
in $\delta$, of the form
\emph{f\textsubscript{\uline{c}}}(\uline{\emph{n}\textsubscript{1}},
\ldots, \emph{\uline{n\textsubscript{k}}}), is replaced with
\emph{f\textsubscript{\uline{b}}}(\emph{f\textsubscript{\uline{c}}}\textsubscript{\uline{1}}(\uline{\emph{n}\textsubscript{1}},
\ldots, \emph{\uline{n\textsubscript{k}}}), \ldots,
\emph{f\textsubscript{\uline{cl}}} (\uline{\emph{n}\textsubscript{1}},
\ldots, \emph{\uline{n\textsubscript{k}}}))). The term that is replaced
in $\delta$ has 1 occurrence of `\emph{f}', whereas the replacing term has
1+\emph{l} occurrences of `\emph{f}'. So the replacement means that the
right side of the newest line has \emph{l} more occurrences of
`\emph{f}' than in $\delta$. Similarly, whereas the counter was set to the
number of `\emph{f}'s in $\delta$, the algorithm now increases it by \emph{l}.
Thus, the counter matches the number of `\emph{f}'s to the right of `='
in the newest line. So, \emph{Case 4} verifies \uline{Lemma 3}.

\emph{Case 5}: The \emph{m}+1\textsuperscript{th} iteration of the Loop applies (R1). Then, the algorithm adds a line
which is the same as the previous one, except that the rightmost, innermost nn-term
in $\delta$, which is of the form
\emph{f\textsubscript{\uline{c}}}(\uline{\emph{n}\textsubscript{1}},
\ldots, \emph{\uline{n\textsubscript{k}}}, \uline{0}), is replaced with
\emph{f\textsubscript{\uline{a}}}(\uline{\emph{n}\textsubscript{1}},
\ldots, \emph{\uline{n\textsubscript{k}}}). The term that is replaced in
$\delta$ has 1 occurrence of `\emph{f}', and so does the term replacing it.
Similarly, the number on the counter is unchanged. Thus, the counter
matches the number of `\emph{f}'s to the right of `=' on the newest
line. So, \emph{Case 5} verifies \uline{Lemma 3}.

\emph{Case 6}: The \emph{m}+1\textsuperscript{th} iteration of the Loop applies (R2). Then, the algorithm adds a line
which is the same as the previous one, except that the rightmost, innermost nn-term
in $\delta$, which is of the form
\emph{f\textsubscript{\uline{c}}}(\uline{\emph{n}\textsubscript{1}},
\ldots, \emph{\uline{n\textsubscript{k}}}, \emph{\uline{m}}) where
\emph{m}\textgreater0, is replaced with the term
\emph{f\textsubscript{\uline{d}}}(\emph{f\textsubscript{\uline{c}}}(\uline{\emph{n}\textsubscript{1}},
\ldots, \emph{\uline{n\textsubscript{k}}}, \uline{\emph{m}-1}),
\uline{\emph{n}\textsubscript{1}}, \ldots,
\emph{\uline{n\textsubscript{k}}}, \uline{\emph{m}--1}). Now the term
that is replaced in $\delta$ has 1 occurrence of `\emph{f}', whereas the term
replacing it has 2. So the replacement means that the right side of the
newest line has one additional occurrence of `\emph{f}'. Similarly,
whereas the counter was set to the number of `\emph{f}'s in $\delta$, the
algorithm now increases it by 1. Thus, the counter matches the number of
`\emph{f}'s to the right of `=' in the newest line. So, \emph{Case 6}
verifies \uline{Lemma 3}.

By induction, this suffices to establish \uline{Lemma 3}. We pause to
note a consequence:
\begin{quote}

\uline{Corollary}: The algorithm halts iff it produces a line
\emph{f\textsubscript{\uline{i}}}(\emph{\uline{i}},
\emph{\uline{n}})=\emph{\uline{m}}, for some \emph{\uline{m}}.

\end{quote}
Proof: The `if' part is true by stipulation, and if \emph{i} is simple,
then the `only if' part is obvious. If \emph{i} is complex, the `only
if' clause is shown as follows: It can be readily verified that the algorithm
halts with complex indices iff the counter reaches 0. From
\uline{Lemma 1}, if \emph{i} is complex, then the algorithm generates
a line \emph{f\textsubscript{\uline{i}}}(\emph{\uline{i}},
\emph{\uline{n}})=$\beta$. And from \uline{Lemma 3}, if it derives
\emph{f\textsubscript{\uline{i}}}(\emph{\uline{i}}, \emph{\uline{n}})=$\beta$,
the counter reaches 0 only if in the latest line of the proof, the term to the right of `=' has zero `\emph{f}'s, i.e., only if it is a numeral. So with complex indices as
well, \uline{Corollary} holds.

The fourth lemma is now stated as follows:
\begin{quote}
\uline{Lemma 4}: If \emph{i} is a complex index, then for any
\emph{n}, if the algorithm derives \emph{f\textsubscript{\uline{i}}}(\emph{\uline{i}}, \emph{\uline{n}})=$\beta$, where $\beta$ is a baseline term, it will continue until it derives \emph{f\textsubscript{\uline{i}}}(\emph{\uline{i}},
\emph{\uline{n}})=\emph{\uline{m}}, for some \emph{\uline{m}}, after which the algorithm halts.
\end{quote}

Proof: Given any complex \emph{i} and \emph{n}, assume the algorithm derives \emph{f\textsubscript{\uline{i}}}(\emph{\uline{i}}, \emph{\uline{n}})=$\beta$, where $\beta$ is a baseline term. As we saw, a baseline term has at least one occurrence of
`\emph{f}'---so at this point, the algorithm will iterate the Main Loop at least
once. Now the Main Loop is defined in such a way that, on a given
iteration, it looks at the equation $\epsilon$ (on the latest line) and to the
right of `=', it finds the rightmost, innermost nn-term $\tau$, if any. If it does not
find such a $\tau$, then the right-hand term of the equation is a numeral and
the algorithm halts---in which case, \uline{Lemma 4} is verified.
Otherwise, if finds such a $\tau$, it adds a line, whereby $\tau$ is replaced in $\epsilon$
with a term that is a p.r. reduction of $\tau$. The counter is then adjusted,
and the algorithm then checks whether the counter is at zero. If so,
then the algorithm halts and, by \uline{Corollary}, we know that the
term right of `=' on the latest line is a numeral. In which case,
\uline{Lemma 4} is verified. Otherwise, Corollary tells us that the Main
Loop will be restarted to further reduce the right-hand term on the
latest line.

So in brief, when the algorithm locates a rightmost, innermost nn-term (to the
right of `='), it adds a line where that term is removed and replaced
with a p.r. reduction of the term. The replacement may not be a numeral,
but no matter: By \uline{Corollary}, the algorithm will continue
iterating and the replacements are made in a unique order, until a line
is reached where the right-hand term is a numeral. And there will be such
a line, given that each p.r. term is ultimately reducible to a numeral
in finitely many steps, as per the definitional rules for p.r.
functions. So the algorithm will produce a derivation from
\emph{f\textsubscript{\uline{i}}}(\emph{\uline{i}},
\emph{\uline{n}})=$\beta$~to
\emph{f\textsubscript{\uline{i}}}(\emph{\uline{i}},
\emph{\uline{n}})=\emph{\uline{m}}, for some \emph{\uline{m}}, by
iterating the Main Loop a sufficient number of times, after which it
will stop---which is what \uline{Lemma 4} says.

From \uline{Lemma 1}, \uline{Lemma 2}, and \uline{Lemma 4}, we know
that if \emph{i} is a complex index, the algorithm yields a unique
\emph{proof} of \emph{f\textsubscript{\uline{i}}}(\emph{\uline{i}},
\emph{\uline{n}})=\emph{\uline{m}}, for some \emph{\uline{m}}. So \textbf{Theorem} holds in the case of complex indices, which
completes the proof of \textbf{Theorem}.

\bigskip
\textbf{6. Why the algorithm is p.r.}

We have just shown that the algorithm behaves as advertised, and we now
present the argument that the algorithm is p.r. 

The algorithm is patently p.r. apart from from steps 4 -- 6 (a.k.a.
``the Main Loop''). When it comes to step 4, the algorithm first applies
Gödel's ``length'' operation (cf. note 8 above) to determine the number
\emph{w} of symbols in the latest line of the proof. This number
\emph{w} is set as a bound for a count of the number of `\emph{f}'s
occurring to the right of `=' in the latest line, thus making the count
of `\emph{f}'s p.r. Moreover, the count of `\emph{f}'s just is the
initial value for the counter, meaning that the counter is initialized
in a primitive recursive manner.

As for step 5, the algorithm first evaluates if the counter's value is
0, in which case it halts. Otherwise, the algorithm can find the
rightmost, innermost nn-term in the latest line of the proof by scanning the term's parentheses structure in a single left-to-right pass to determine nesting depth. This is a bounded operation on a finite string and appears primitive recursive. In more detail, the algorithm counts the parenthesis-pairs---not greater than \emph{w}---which
enclose each `\emph{f}' to the right of `='. An occurrence of `\emph{f}'
that yields the greatest count begins an innermost nn-term.\footnote{Note
  that codes for terms in PRA-- can be identified by Gödel's
  ``first-order term'' operation; \#18 in Gödel's list of p.r. functions
  and relations.} If there is more than one, the rightmost one is the last one in the wff. (If the innermost terms in the wff are enumerated
left-to-right, it will be the one enumerated by the greatest
number, not greater than \emph{w}.) Again, the relevant maneuvers are p.r.;
they amount to bounded arithmetic operations on finite sequences.\footnote{Assistance with this paragraph was provided by GPT 4.0.}

At step 6, the algorithm adds a new line to the proof, which is
identical to the previous line, except that the rightmost, innermost
nn-term $\tau$ has been replaced. Importantly, replacement is p.r.; cf.
Gödel's \emph{Sb} operation, \#31 in his list. In more detail, the algorithm first identifies the index for $\tau$ as 0, 1, or as coding a specific sort of string. (If the last, the first element of the string indicates if the index is for projection function, or if it is a composition or a recursion index.) From this determination, the algorithm decides on one of the replacements listed below. Since there are only finitely many replacement templates, this decision process may be implemented by a p.r. predicate. 

\begin{quote}
Per (Z): If $\tau$ is \emph{f}\textsubscript{\uline{0}}(\emph{\uline{m}}),
replace it with \uline{0}.

Per (S): If $\tau$ is \emph{f}\textsubscript{\uline{1}}(\emph{\uline{m}}),
replace it with \emph{\uline{m}}$^\prime$.

Per (P): If $\tau$ is
\emph{f\textsubscript{\uline{c}}}(\uline{\emph{n}\textsubscript{1}},
~\ldots{}, \emph{\uline{n\textsubscript{k}}}) where \emph{c} codes
\emph{\uline{j}}*\emph{\uline{k}} and $1\!\leq\!\emph{j}\!\leq\!\emph{k}\!\leq\!\emph{c}$,
replace it with \emph{\uline{n\textsubscript{j}}}.

Per (C): If $\tau$ is
\emph{f\textsubscript{\uline{c}}}(\uline{\emph{n}\textsubscript{1}},~
~\ldots{}, \emph{\uline{n\textsubscript{k}}}), and \emph{c} is a
composition index, replace it with
\emph{f\textsubscript{\uline{b}}}(\emph{f\textsubscript{\uline{c}}}\textsubscript{\uline{1}}(\uline{\emph{n}\textsubscript{1}},
\ldots, \emph{\uline{n\textsubscript{k}}}), \ldots,
\emph{f\textsubscript{\uline{cl}}}(\uline{\emph{n}\textsubscript{1}},
\ldots, \emph{\uline{n\textsubscript{k}}})).

Per (R1): If $\tau$ is
\emph{f\textsubscript{\uline{c}}}(\uline{\emph{n}\textsubscript{1}},
~\ldots{}, \emph{\uline{n\textsubscript{k}}}, \uline{0}), and \emph{c}
is a recursion index, replace it with
\emph{f\textsubscript{\uline{a}}}(\uline{\emph{n}\textsubscript{1}},
~\ldots{}, \emph{\uline{n\textsubscript{k}}}).

Per (R2): If $\tau$ is
\emph{f\textsubscript{\uline{c}}}(\uline{\emph{n}\textsubscript{1}},
~\ldots{}, \emph{\uline{n\textsubscript{k}}}, \uline{\emph{m}}$^\prime$), and
\emph{c} is a recursion index, replace it with
\emph{f\textsubscript{\uline{d}}}(\emph{f\textsubscript{\uline{c}}}(\uline{\emph{n}\textsubscript{1}},
\ldots, \emph{\uline{n\textsubscript{k}}}, \uline{\emph{m}}),
\uline{\emph{n}\textsubscript{1}}, \ldots,
\emph{\uline{n\textsubscript{k}}}, \uline{\emph{m}}).

\end{quote}

In each case, the replacement operation is bounded by the length of the term, and the output string can be constructed explicitly using only primitive recursive operations. In detail, given \emph{c}, the algorithm extracts \emph{\uline{j}},
\emph{\uline{k}}, \emph{\uline{b}}, \uline{\emph{c}\textsubscript{1}},
~\ldots{}, \emph{\uline{c\textsubscript{l}}}, \emph{\uline{a}}, or
\emph{\uline{d}}, as needed, and such extraction is p.r. And the replacing term
for $\tau$ is constructed in a p.r. manner since this at most requires
replacement into finite term-schemes (which are found in the
axiom-schemes for the function symbols). 

After replacing $\tau$, the counter is adjusted (or not) based on the
rule applied. For (Z), (S), and (P), it is decreased by 1, and for (R1)
and (R2), respectively, it is unchanged or increased by 1. As for (C),
it is increased by \emph{l}, the number of sub-functions in the
composition (which is determined by a count bounded by its index). These updates are p.r.; at most they require addition or
subtraction on a pair of numbers. And the counter eventually reaches 0; cf. \uline{Lemma 3}.

So too, the other operations in the algorithm are patently p.r. (e.g.,
halting when the counter reaches zero). It thus appears that the procedure never searches indefinitely using unbounded minimization. The proof checking is implemented by direct pattern-matching against a fixed finite set of templates, and the procedure always operates over terms of bounded length.

In short, while the algorithm computes a function that would ordinarily be excluded by standard diagonal arguments, the procedure appears to conform to the criteria for primitive recursive definition.

\bigskip
\textbf{7. Closing}

The algorithm presented in this paper appears to satisfy all standard criteria for being primitive recursive while enabling the behavior of a universal p.r. function. This challenges the standard view that diagonalization tells against the existence of such a function. Other diagonal arguments are known to demonstrate ill-defined sets rather than the impossibility of a diagonal object (Simmons, op. cit.), and the present construction indicates that diagonalization in the present case is of this sort. Still, this leaves open why the diagonal argument might not foreclose the possibility of a universal p.r. function. Clarifying this could provide deeper insight into primitive recursive definitions or diagonalization within formal arithmetic. Either way, the case suggests that further analysis may refine our understanding of what diagonal arguments presuppose about the structure of definability or primitive recursion.\footnote{My thanks to Bill Gasarch, GPT 4.0, Bill Mitchell, Panu Raatikainen, Lionel Shapiro, Henry Towsner, Nic Tideman, Bruno Whittle, Noson Yanofsky, and Richard Zach for discussion of issues relevant to this paper. I also express gratitude to an audience at the 2022 meeting of the Australasian Association of Philosophy.}

\bigskip
\textbf{Works Cited}

Curry, H. (1941). `A Formalization of Recursive Arithmetic,'
\emph{American Journal of Mathematics} 63: 263-282.

Gödel, K. (1931). `Über Formal Unentscheidbare Sätze der \emph{Principia
Mathematica} und Verwandter Systeme I,'~\emph{Monatshefte für Mathematik
Physik} 38: 173--198. Pagination is from Gödel, K. (1986).
\emph{Collected Works I. Publications 1929--1936}. S. Feferman et al.
(eds.), Oxford: Oxford University Press, pp. 144--195.

Hilbert, D. \& Bernays, P. (1934).~\emph{Grundlagen der Mathematik},
\emph{Vol. I}, Berlin: Springer.

Kleene, S.C. (1952). \emph{Introduction to Metamathematics}. Amsterdam:
North Holland Publishing Co.

Nelson, E. (2011). \emph{Elements}.
\url{https://arxiv.org/abs/1510.00369}

Nelson, E. (2011a). `Re: The Inconsistency of Arithmetic.' \emph{The
n-Category Café}.
\url{https://golem.ph.utexas.edu/category/2011/09/the_inconsistency_of_arithmeti.html\#c039590}.

Quine, W. (1951). \emph{Mathematical Logic}, revised edition. Cambridge,
MA: Harvard University Press.

Robinson, R.M. (1950). `An Essentially Undecidable Axiom System,'
\emph{Proceedings of the International Congress of Mathematics}:
729-730.

Simmons, K. (1993). \emph{Universality and the Liar}. Cambridge: Cambridge University Press.

Skolem, T. (1923). `Begründung Der Elementaren Arithmetik Durch Die
Rekurrierende Denkweise Ohne Anwendung Scheinbarer Veranderlichen Mit
Unendlichem Ausdehnungsbereich',~\emph{Videnskapsselskapets Skrifter, I.
Matematisk-Naturvidenskabelig Klasse}, 6: 1--38. Translated by S.
Bauer-Mengelberg as `The Foundations of Elementary Arithmetic
Established by the Recursive Mode of Thought, without the Use of
Apparent Variables Ranging over Infinite Domains.' In Heijenoort, J. van
(ed.). (1967), \emph{From Frege to Gödel: A Source Book in Mathematical}
\emph{Logic, 1879--1931}. Cambridge, MA: Harvard University Press, pp.
302-333.

Tao, T. (2011) `Re: The Inconsistency of Arithmetic.' \emph{The
n-Category Café}.
\url{https://golem.ph.utexas.edu/category/2011/09/the_inconsistency_of_arithmeti.html\#c039553}.

\end{document}